\documentclass[a4paper,10pt]{article}

\usepackage{tikz}
\usepackage{pgfplots}
\pgfplotsset{compat=1.14}

\usepackage{a4wide}
\usepackage[utf8]{inputenc}
\usepackage[T1]{fontenc}
\usepackage{amsmath,amssymb,commath,mathtools,amsfonts}
\usepackage{mathrsfs}
\usepackage{enumitem}
\usepackage{booktabs}
\usepackage{subcaption}
\usepackage{siunitx}

\usepackage[pdfusetitle]{hyperref}

 \hypersetup{
   colorlinks=true,
  citecolor=blue,
  linktocpage=true,
  pdfstartview=FitV,
  breaklinks=true,
  pdfpagemode=UseNone,
  pageanchor=true,
  pdfpagemode=UseOutlines,
  plainpages=false,
  bookmarksnumbered,
  bookmarksopen=true,
  bookmarksopenlevel=1,
  hypertexnames=false,
  pdfhighlight=/O,
  nesting=true,
  pdfsubject={Mathematics, Numerical analysis},
  pdfkeywords={power networks, Lyapunov function, stability, numerical approximations},
  pdfcreator={pdfLaTeX},
  pdfproducer={},
}

\usepackage{amsthm,thmtools}
\declaretheorem[name=Theorem]{thm}
\declaretheorem[name=Proposition, sibling=thm]{prop}
\declaretheorem[name=Remark, sibling=thm]{remark}
\declaretheorem[name=Lemma, sibling=thm]{lemma}

\usepackage[backend=bibtex,style=alphabetic,sortcites=true,giveninits=true]{biblatex}  %
\DeclareMathOperator{\diag}{diag}
\DeclareMathOperator{\supp}{supp}

\DeclareMathOperator{\range}{range}
\DeclareMathOperator{\ODE}{ODE}
\DeclareMathOperator{\PDE}{PDE}
\DeclareMathOperator{\TV}{TV}
\DeclareMathOperator{\sign}{sign}

\let\abs\relax
\DeclarePairedDelimiter\abs{\lvert}{\rvert}

\title{On the relation of powerflow and Telegrapher's equations: continuous and numerical Lyapunov stability}
\author{E. Fokken\footnotemark[1] \and S. G\"ottlich\footnotemark[1]}

\begin{document}

\maketitle
\renewcommand*{\thefootnote}{\fnsymbol{footnote}}
\footnotetext[1]{Department of Mathematics, University of Mannheim, \{\href{mailto:fokken@uni-mannheim.de}{fokken},\href{mailto:goettlich@uni-mannheim.de}{goettlich}\}@uni-mannheim.de}
\renewcommand*{\thefootnote}{\arabic{footnote}}%

\abstract{In this contribution we analyze the exponential stability of power networks modeled with the Telegrapher's equations as a system of balance laws on the edges. We show the equivalence of periodic solutions of these Telegrapher's equations and solutions to the well-established powerflow equations.
In addition we provide a second-order accurate numerical scheme to integrate the powerflow equations and show (up to the boundary conditions) Lyapunov stability of the scheme.}\\

[0.5ex] %

\section{Introduction}
\label{sec:introduction}

In recent years, due to the need to restructure energy systems to incorporate more renewable energy sources renewed interest into energy-systems has sparked.
In the case of electric transmission lines several mathematical approaches can be found in the literature, see for example \cite{Andersson15book,doi:10.1137/1.9781611974164,Frank16,Gottlich2021} for an overview. 
All modeling approaches rely on a graph structure where generators of electrical powers and consumers at nodes are connected by transmission lines.
From a physical point of view, voltage and current are transported along the lines while power loss might occur due to resistances.
In many applications, the so-called power flow equations provide a well-established  tool to analyze the performance of electric transmission lines over time.
Mathematically, the resulting nonlinear system of equations is typically solved via Newton's method.
Another approach to study not only the temporal resolution of power flow but also the spatial resolution in transmission lines are the spatially one-dimensional Telegrapher's equations.
These equations are a coupled system of linear hyperbolic balance laws for voltage and current.
In this work, we aim to associate these two modeling approaches and therefore rigorously analyze the relation between the powerflow equations used to model power systems and the Telegrapher's equations, that describe the underlying physics in more detail.

In terms of operation of electric transmission lines the investigation of stability issues is of great importance.
Recent research in this direction include exponential stability considerations \cite{Egger_2017} as well as open- and closed loop problems.
Open loop problems tackle questions of optimal electricity input \cite{Gottlich2019,GottlichPT2019, Gottlich2018},  while closed loop control are concerned with feedback control \cite{Gugat2014} or boundary stabilization \cite{Gugat_2018}.
Similar to the ideas presented in \cite{Gottlich2016}, we study the stability of the Telegrapher's equations using the concept of Lyapunov functions \cite{bastin2016stability}.
We prove a stability result for the Telegrapher's equations on networks with (linear) coupling conditions and show how such a solution can be thought of as a solution to the powerflow equations.
In particular, we are interested in how a numerical approximation might mimic the Lyapunov stability as for example introduced in \cite{bastin2016stability}.
In contrast to the already existing literature \cite{Gerster2019,Gottlich2016,goettlich_schillen_2017}, where only first order schemes have been applied, we introduce a numerical scheme of second order to examine the stability. A numerical analysis shows that the proposed scheme is Lyapunov stable and can be used to
numerically study power networks as the simulation results indicate.

The paper is organized as follows: In section \ref{sec:transm-line-model} we present the physical meaning of the Telegrapher's equation, define the network setting and prove theoretical properties concerning the Lyapunov stability. Section \ref{sec:power-network} is concerned with the derivation of the well-known powerflow equations from the Telegrapher's equation with a particular focus on the power coupling conditions. In section \ref{sec:numercial-scheme}, a numerical scheme of second order is introduced to solve the Telegrapher's equation. The Lyapunov stability for the numerical approximation is also discussed in detail. The last section \ref{sec:numerical-examples} deals with numerical experiments to investigate the performance of the scheme.

\section{Solution structure of the Telegrapher's equations}%
\label{sec:transm-line-model}

The Telegrapher's equation is a system of linear hyperbolic partial differential equations (PDE) with linear source term in one space dimension and therefore often called {\em balance law}, see \cite{Gottlich2016,Jiang2005} for more details.
It models the time-dependent evolution of voltage and current on a transmission line and can be written in terms of an initial-boundary value problem on a single line with $x\in[0,\ell]$ and $t\in[0,T]$ as:
\begin{equation}
  \label{eq:6}
  \begin{aligned}
    &u_t +Au_x +Fu=0,\\
    &u(0,x) = u_0(x),\\
    &u^1(t,0)=v_0(t),\\
    &u^1(t,\ell)=v_\ell(t),
\end{aligned}
\end{equation}
where $u:[0,T]\times[0,\ell] \to \mathbb{R}^2$ and
\begin{equation}
  \label{eq:4}
  A=
  \begin{pmatrix}
    0&\frac{1}{C}\\
    \frac{1}{L}&0
  \end{pmatrix},\
  F=
  \begin{pmatrix}
    \frac{G}{C}&0\\
    0&\frac{R}{L}
  \end{pmatrix}.
\end{equation}
Here, $R,L,G,C>0$ are constants depending on the line, $T>0$ is some time horizon and $\ell$ is the length of the line.
The vector $u$ is composed of the voltage $u^1=v$ and the current $u^2=i$.
In the engineering literature (e.g.~\cite[§1.2.2 Lossy Transmission Lines]{transmissionlinesandlumpedcircuits}) the Telegrapher's equation is usually written as
\begin{equation*}
  \begin{aligned}
    \pd{v}{x} &= -(Ri + L\pd{i}{t}),\\
    \pd{i}{x} &= -(Gv + C\pd{v}{t}).
  \end{aligned}
\end{equation*}

Yet another form of the Telegraper's equation is gotten by transformation into characteristic variables $\xi^+,\xi^-$:
\begin{equation}
  \label{eq:14}
  \begin{aligned}
    \sqrt{\tfrac{C}{L}}v &=\xi^+-\xi^-,\\
    i &=\xi^++\xi^-,
  \end{aligned}
\end{equation}
so that equation~\eqref{eq:6} turns into
\begin{equation*}
  \xi_t +\Lambda \xi_x +B\xi = 0,
\end{equation*}
with
\begin{equation}
  \label{eq:1}
  \Lambda =\frac{1}{\sqrt{LC}}
  \begin{pmatrix}
    1& \\ &-1
  \end{pmatrix}, \
  B =
  \begin{pmatrix}
    a & b\\
    b & a
  \end{pmatrix}
\end{equation}
and $a = \frac{R}{L}+\frac{G}{C}$, $b=\frac{R}{L}-\frac{G}{C}$.
Note that the eigenvalues of $B$ are $a+b>0$ and $a-b>0$ and that in addition $a>0$ holds.

In a power network setting, these equations are defined on the lines and are coupled at the nodes by some conditions.
The overall aim of this work is to show that in so-called alternating current (\emph{AC}) networks ---
which are almost all power networks in use in energy systems --- the solution normally used in engineering, namely the solution to the powerflow equations, can be thought of as an exponentially stable solution of the Telegrapher's equations defined on the power network.
In an AC power network all voltages and currents that serve as boundary conditions are sinusoidal with the same angular frequency $\omega$, that is, they are of the form $A \sin(\omega t) + B \cos(\omega t)$ for real $A,B$ and $\omega$.

We start of by examining a time-periodic solution of the Telegrapher's equations on a line with the same angular frequency $\omega$ as that of the boundary conditions.
Therefore we write $v$ and $i$ as Fourier series in time,
\begin{equation}
  \label{eq:2}
  \begin{aligned}
    v(x,t) &= \sum_{m \in \mathbb{Z}}V_m(x) e^{jm \omega t},\\
    i(x,t) &= \sum_{m \in \mathbb{Z}}I_m(x) e^{jm \omega t},\\
  \end{aligned}
\end{equation}
with complex-valued functions $V_m,I_m$, which we call complex voltage and complex current.
For the voltage and current to be real-valued, we impose
\begin{equation*}
  V_m =   V_{-m}^*\text{ and }I_m =   I_{-m}^*\text{ for all }m \in \mathbb{Z}.
\end{equation*}

Due to the linear independence of the system ${(e^c)}_{c \in \mathbb{C}}$ the Telegrapher's equations decompose into a decoupled family of complex linear ordinary differential equations (ODEs), indexed by $m \in \mathbb{Z}$, which only governs $V_m$ and $I_m$ via
\begin{equation}
  \label{eq:7}
  \begin{aligned}
    \pd{V_m}{x} &= -(R + jm\omega L)I_m,\\
    \pd{I_m}{x} &= -(G + jm\omega C)V_m .
  \end{aligned}
\end{equation}
We now define
\begin{equation*}
  \begin{aligned}
    Y^0_m &= \frac{\sqrt{G + j\,m\,\omega\,C}}{\sqrt{R + j\,m\,\omega\,L}},\\
    \gamma_m &= \sqrt{(R + j\,m\,\omega\,L)}\cdot \sqrt{(G + j\,m\,\omega\,C)},
  \end{aligned}
\end{equation*}
and it is easily verified, that
\begin{equation}
  \label{eq:10}
  \begin{aligned}
    {V_m}(x) &= \phantom{-}\frac{1}{\sinh(\gamma_m \ell)}\left[ {V}^\ell_m \sinh(\gamma_m x) +  {V}^0_m \sinh(\gamma_m (\ell-x)) \right],\\
    {I_m}(x) &= -\frac{Y^0_m}{{\sinh(\gamma_m \ell)}} \left[  {V}^\ell_m \cosh(\gamma_m x) - {V}^0_m \cosh(\gamma_m (\ell-x)) \right]
  \end{aligned}
\end{equation}
solve equations~\eqref{eq:7}.  Here $V^0_m, V^\ell_m$ are the complex voltages at either sides of the line, which we choose as boundary conditions for now.
It is also possible to prescribe the current at the end of a line as done in the network setting in \autoref{sec:stab-telegr-equat}.

Note also that $Y^0_m$ and $ \sinh(\gamma_m \ell)$ are never zero because $R,L,G,C>0$.
Plugging these into equation~\eqref{eq:2} yields a solution to the Telegrapher's equations. This solution will be the basis for deriving the powerflow equations in Section~\ref{sec:power-network}. Next, we will investigate the stability of the Telegrapher's equation on networks.

\subsection{Stability of the Telegrapher's equations on networks}%
\label{sec:stab-telegr-equat}

Before we can study the stability of the Telegrapher's equations for power networks, we have to introduce the network setting and coupling conditions.
Stability of similar networks has been analyzed,
see e.g.~\cite{Egger_2017,Gugat_2018,Gottlich2016,goettlich_schillen_2017}.
Like the authors of~\cite{Egger_2017,Gugat_2018}, we also also apply the formalism of operator semigroups to our network to show well-posedness of our problem.
Similar to \cite{Egger_2017} we show the stability of the system without feedback control. 

We discuss two different stability estimates depending on the network.
The first one yields $L^2$-stability for the networks we examine.
The second one yields stability in $\norm{\cdot}_\infty$ but is valid only for a network with equal line parameters. 

In the end we want to model power networks, that is, generators of electrical powers and consumers (loads), that are connected by transmission lines.
Hence we use a directed graph as a power network model.
Let therefore $\mathscr{G} = (\mathscr{V},\mathscr{E})$ be a finite directed graph with vertex set $\mathscr{V}$ and edge set $\mathscr{E}$.
An edge $e$ consists of an interval $[0,\ell_e]$ together with edge parameters $R_e,L_e,G_e,C_e$.
To a vertex we assign a virtual point $x_v$ and
for an edge $s$ starting at $v$ we identify $0 \in [0,\ell_s]$ with $x_v$.
For an edge $e$ ending in $v$ we identify $\ell_s \in [0,\ell_s]$ with $x_v$.
The reasoning behind this identification is to interpret the graph as a subset of some $\mathbb{R}^n$, where the start and endpoints of an edge
are actually the points where the vertices are located.
This image is most natural for planar graphs, where $\mathbb{R}^2$ can be used.

Let further $\mathscr{V} = \mathcal{G}\ \dot{\cup}\ \mathcal{L}$ where $\mathcal{G}$ denotes the nodes with a generator and $\mathcal{L}$ denotes the
nodes with a load so that every node is either a generator or a load.
Also for $v \in \mathscr{V}$ let $\mathscr{E}_v\subset \mathscr{E}$ be the set of edges connected to $v$.
In addition we define the function
\begin{equation}
  \label{eq:81}
  s(e,v) = \begin{cases}
    1 \text{ if $v$ is the end of $e$},\\
    -1 \text{ if $v$ is the start of $e$},
  \end{cases}
\end{equation}
which distinguishes start and end nodes of an edge.

Then we formulate the linear inhomogeneous initial-boundary-value problem~\eqref{IIBVP} on the graph $\mathscr{G}$,
\begin{equation}
  \tag{IIBVP}\label{IIBVP}
  \partial_t u_e = -A_e\partial_{x} u_e -F_e u_e\ \forall e \in \mathscr{E},
\end{equation}
with matrices $A_e$ and $F_e$ of the structure of~\eqref{eq:4} and
under the conditions
\begin{align}
  u_e(0,x) &={(u_0)}_e(x) \forall e \in \mathscr{E} \tag{IC}\label{IC},\\
  u^1_e(t,x_v) &= u^1_f(t,x_v)\ \forall\ e,f: x_v \in e \tag{CC}\label{CC},\\
  u^1_g(t,x_g) &= v_g(t)\ \forall\ g \in \mathcal{G}\ \forall \ e \in \mathscr{E}_g \tag{KCg}\label{KCg},\\
  \sum_{e \in \mathscr{E}_l}s(e,l)u^2_e(t,x_l) &= i_l\ \forall\ l \in \mathcal{L} \tag{KCl},\label{KCl}
\end{align}
where~\eqref{IC} is the initial condition,~\eqref{CC} is a continuity condition for the voltage in the nodes,
\eqref{KCg} is a coupling condition fixing the voltage at generators and~\eqref{KCl} is a coupling condition fixing the current at loads.
In order to show well-posedness we use the theory of operator semigroups, see~\cite{engel_one_parameter_semigroups}
and especially~\cite{semigroups_boundary}.
Therefore we introduce the following setting:
As underlying Banach (and actually Hilbert) space we take
\begin{equation*}
  X = L^2(\prod_{e \in \mathscr{E}} [0,\ell_e],\mathbb{R}^2),
\end{equation*}
with the scalar product
\begin{equation}
  \label{eq:79}
  < u, v > = \frac{1}{2}\sum_{e \in \mathscr{E}}\int_0^{\ell_e} (C_eu_e^1v_e^1 + L_eu^2_ev^2_e)\dif x .
\end{equation}
As linear operator we take
\begin{equation*}
  \begin{aligned}
    \mathcal{A}:D(\mathcal{A})&\to X\\
    u &\mapsto \left(e \ni x \mapsto -A_e\partial_x u_e(x)- F_eu_e\right)
  \end{aligned}
\end{equation*}
and as its domain the space
\begin{equation}
  \label{eq:78}
  D(\mathcal{A}) = H^1(\prod_{e \in \mathscr{E}} [0,\ell_e],\mathbb{R}^2).
\end{equation}
In addition we choose the boundary operator
\begin{equation*}
  \begin{aligned}
    \mathscr{L}: D(\mathcal{A})&\to \mathbb{R}^{2\abs*{\mathscr{E}}},\\
    u&\mapsto (B_V,B_C,B_I),
  \end{aligned}
\end{equation*}
with
\begin{equation*}
  \begin{aligned}
    B_V &= (u^1_e(t,x_g))_{g \in \mathcal{G}, e \in \mathscr{E}_g},\\
    B_C &=  ((u^1_e(x_l) - u^1_{e+1}(x_l))_{e = 1}^{\abs{\mathscr{E}_l}-1})_{l \in \mathcal{L}},\\
    B_I &= \left(\sum_{e \in \mathscr{E}_l}s(e,l)u^2_e(t,x_l)\right)_{l \in \mathcal{L}},
  \end{aligned}
\end{equation*}
which maps a collection of $H^1$-functions on the edges to the evaluation of line voltages at generator nodes, the differences in line voltages at load nodes and load currents, respectively.
Note that $\mathscr{L}$ is continuous on $D(\mathcal{A})$ as a concatenation of the (continuous) trace operator and a finite-dimensional linear mapping.

To prove well-posedness we first verify the prerequisites of \cite[Assumption 3.1]{semigroups_boundary} in
\begin{lemma}
  \label{lemma:assumptions}
  The operator $\mathcal{A}$ satisfies
  \begin{enumerate}[label=(G\arabic*), ref=(G\arabic*)]
  \item \label{item:G1}The restriction of $\mathcal{A}$ to $\ker \mathscr{L}$, namely $\mathcal{A}_0 = \mathcal{A}|_{\ker \mathscr{L}}: D(\mathcal{A}_0) = \ker \mathscr{L} \to X$ is densely defined, closed and has non-empty resolvent set.
  \item \label{item:G2}$\mathscr{L}$ is surjective.
  \item \label{item:G3}The combined operator
    \begin{equation*}
      \begin{pmatrix}
        \mathcal{A}\\
        \mathscr{L}
      \end{pmatrix}: D(\mathcal{A}) \to X \times \mathbb{R}^{2\abs*{\mathscr{E}}}
    \end{equation*}
    is closed.
  \end{enumerate}
  \begin{proof}
    For this proof we consider the basis introduced in~\eqref{eq:14}, where $\mathcal{A}$ is of the form
    \begin{equation*}
          \mathcal{A} = -\lambda
      \begin{pmatrix}
        1\\&-1
      \end{pmatrix}\partial_x -
      \begin{pmatrix}
        a&b\\b&a
      \end{pmatrix}.
    \end{equation*}
    For \ref{item:G1} we first note that the space with zero boundary conditions, i.e.,
    \begin{equation*}
      H^1_0(\prod_{e \in \mathscr{E}} [0,\ell_e],\mathbb{R}^2),
    \end{equation*}
    is a subset of $D(\mathcal{A})\cap \ker \mathscr{L}$, which is well-known to be a dense subspace of $L^2(\prod_{e \in\ \mathscr{E}} [0,\ell_e],\mathbb{R}^2)$.
    Then we note that $\mathcal{A}$ is just the sum of the spatial derivative and a bounded linear operator.
    As $\mathscr{L}$ is continuous on $H^1$, its kernel is a closed subspace therein.
    Together we find that $\mathcal{A}_0$ is closed, if the derivative operator $\partial_x:H^1 \to L^2$ is closed, which is well-known.
    For the non-empty resolvent set we show that $0$ is in the resolvent set of $\mathcal{A}_0$,
    that is $\mathcal{A}_0$ is injective and its left inverse $\mathscr{R}$ is densely-defined and bounded.
    The domain of $\mathscr{R}$ is just the range of $\mathcal{A}_0$, so we show this range to be dense in $X$.
    Therefore we show that the space of test functions
    $\mathcal{T} = \set{f \in \prod_{e \in \mathscr{E}}\in C^{\infty}([0,\ell_e],\mathbb{R}^2)| \supp(f) \subset [0,\ell_e]}$ fulfills $\mathcal{T} \subset \range(\mathcal{A}_0)$ as it is well-known that $\mathcal{T}$ is dense in $L^2$.

    Let $g \in \mathcal{T}$.  There holds $\mathcal{T} \subset H^1_0$ and as the zero boundary conditions essentially decouple different edges, we only consider one edge, given by $[0,\ell]$.
    The Fourier transform of $g$ (which is well-defined as $g \in L^2$) is given by
    \begin{equation*}
        g(x) = \int_{-\infty}^{\infty}
        \begin{pmatrix}
          \hat{g}_1(k)\\
          \hat{g}_2(k)
        \end{pmatrix}
        e^{ikx} \dif k . 
    \end{equation*}
    The operator $\mathcal{A}_0$ evaluated on $g$ is given by
    \begin{equation*}
      g(x) = -\int_{-\infty}^{\infty}
      \begin{pmatrix}
        a + \lambda ki & b\\
        b &a - \lambda ki
      \end{pmatrix}
      \begin{pmatrix}
        \hat{g}_1(k)\\
        \hat{g}_2(k)
      \end{pmatrix}
      e^{ikx} \dif k.
    \end{equation*}
    The matrix is invertible for all $k$, because $a^2-b^2>0$.
    Let us call this matrix $D$.
    If we choose
    \begin{equation*}
      h(x) = \int_{-\infty}^{\infty}D^{-1}
      \begin{pmatrix}
        \hat{g}_1(k)\\
        \hat{g}_2(k)
      \end{pmatrix}
      e^{ikx} \dif k,
    \end{equation*}
    we have $\mathcal{A}_0h = g$.
    Then we get $h \in L^2([0,\ell])$ and $h \in C^{\infty}([0,\ell])$. It remains to show that $h$ is again supported in $[0,\ell]$.
    This is the case because of the Paley–Wiener theorem, see \cite[Theorem 19.3]{rudin1987real}, which relates $L^2$-functions supported on an interval and their Fourier transforms.
    Together we have $h \in H^1_0([0,\ell]) \subset D(\mathcal{A}_0)$.
    Finally using $\norm{g}_2 = \norm{\hat{g}}_2$ and the definition of $D$, it is easily computed that $\mathscr{R}$ is bounded.
    The injectivity of $\mathcal{A}_0$ will be shown in \autoref{lemma:injectivity}.

    For \ref{item:G2} we note that functions, that on every edge have the form
    \begin{equation*}
      u(x) =
      \begin{pmatrix}
        a+bx\\
        c+dx
      \end{pmatrix},
    \end{equation*}
    are in $D(\mathcal{A})$ and can be used to construct any value for $\mathscr{L}$.

    For \ref{item:G3} we note again that $\mathscr{L}$ is continuous on $D(\mathcal{A})$ and that $\mathcal{A}$ is closed, hence their combination is as well.
  \end{proof}
\end{lemma}

\begin{lemma}
  \label{lemma:injectivity}
  $\mathcal{A}_0$ is injective.
  \begin{proof}
    Here we consider again the basis, where the components of $u$ are the voltage $u^1 = v$ and the current $u^2 = i$. 
    Therefore consider an element $u_0 = (v,i)^T \in \ker \mathcal{A}_0$.
    Such $u_0$ is a valid initial condition for~\eqref{IIBVP} with homogeneous coupling conditions and a corresponding solution given by $u(t,x) = u_0(x)$ for all $t \geq 0$.
    This solution is at least once differentiable with respect to time and space.
    Therefore we can apply the reasoning of the following \autoref{prop:4} to this special solution and find
    \begin{equation*}
      0 = \sum_{e \in \mathscr{E}}\int_0^{\ell_e}C_e v_e\dot{v}_e + L_e i_e\dot{i}_e =     \dot{\mathcal{V}}(\zeta(t)) \leq - 2\min_{e \in \mathscr{E}}\left( \min\left( \frac{R_e}{L_e},\frac{G_e}{C_e} \right) \right) \mathcal{V}(u_0).
    \end{equation*}
    As $\mathcal{V}(u_0) = 0$ if and only if $u_0 = 0$, we find that $\mathcal{A}_0$ is injective.
  \end{proof}
\end{lemma}

\begin{lemma}
  \label{lemma:2}
  The~\eqref{IIBVP} is well-posed and therefore admits a solution $u \in D(\mathcal{A})$ for $u_0 \in D(\mathcal{A})$ and suitably smooth $v_g$, $i_l$.
  \begin{proof}
    Following \cite{semigroups_boundary}, it remains to show that $\mathcal{A}_0$ generates a strongly continuous semigroup.
    For this we first omit the matrices $F_e$ in $\mathcal{A}_0$ and note that the resulting operator is skew-symmetric with respect to the scalar product~\eqref{eq:79} of $X$. (The boundary terms arising in integration by parts vanish, because $D(\mathcal{A}_0) \subset \ker \mathscr{L}$).
    Then we cite \cite[3.24 Theorem. (Stone, 1932)]{engel_one_parameter_semigroups} to see that $\mathcal{A}_0$ generates a unitary group.
    We now add $F_e$ again, which is a bounded linear operator and use \cite[Proposition 1.12]{engel_one_parameter_semigroups},
    showing that the full operator $\mathcal{A}_0$ also generates a strongly continuous semigroup.
    The last step of non-vanishing $v_g,i_l$ follows from application of the variation of constants formula~\cite[VI, Corollary 7.8]{engel_one_parameter_semigroups}, which introduces some smoothness conditions on $v_g,i_l$.
    The exact nature of these condition is of no consequence to us but could be found by working through~\cite{semigroups_boundary} to incorporate the boundary conditions as initial conditions together with a source term.    
  \end{proof}
\end{lemma}

With \autoref{lemma:2} we have a $H^1$-solution to~\eqref{IIBVP}.
We now show its exponential stability using the Lyapunov function
\begin{equation*}
  \mathcal{V}(u(t)) = \frac{1}{2} \sum_{e \in \mathscr{E}} \int_0^{\ell_e}C_e(u^1_e)^2 + L_e(u^2_e)^2 .
\end{equation*}

\begin{prop}\label{prop:4}
  Two solutions $u, w$ of the~\eqref{IIBVP} with identical coupling conditions $v_g,i_l$ but possibly different initial conditions $u_0,w_0$ grow closer exponentially in the $L^2$-norm.
  \begin{proof}
We take the difference of the two solutions $\zeta = u-w = (v, i)^T$ with initial condition $\zeta_0 = u_0 -w_0$ and examine the derivative of the Lyapunov function.
    As before the coupling condition for $\zeta$ is then homogeneous.
    \begin{equation}
      \label{eq:80}
      \begin{aligned}
        \dot{\mathcal{V}}(\zeta(t)) &= \sum_{e \in \mathscr{E}}\int_0^{\ell_e}C_e v_e\dot{v}_e + L_e i_e\dot{i}_e\\
        &= \sum_{v \in \mathscr{V}}\sum_{e \in \mathscr{E}_v}v_e(t,x_v)s(e,v)i_e(t,x_v) - \sum_{e \in \mathscr{E}}\int_0^{\ell_e}G_e v_e^2 + R_e i_e^2\\
        &= \sum_{v \in \mathscr{V}}v_v(t,x_v)\sum_{e \in \mathscr{E}_v}s(e,v)i_e(t,x_v) - \sum_{e \in \mathscr{E}}\int_0^{\ell_e}G_e v_e^2 + R_e i_e^2\\
        &= - \sum_{e \in \mathscr{E}}\int_0^{\ell_e}G_e v_e^2 + R_e i_e^2\\
        &\leq - 2\min_{e \in \mathscr{E}}\left( \min\left( \frac{R_e}{L_e},\frac{G_e}{C_e} \right) \right) \mathcal{V}(\zeta(t)),
      \end{aligned}
    \end{equation}
    where the third equality stems from~\eqref{CC} and the fourth from~\eqref{KCg} or~\eqref{KCl} depending on the node being either a load or a generator.  Note that because of the homogeneous boundary conditions, for generators $v$ there holds $v_v(t,x_v) = 0$ and for loads there holds $\sum s(e,v)i_e(t,x_v) = 0$.
    This, together with the fact that $\mathcal{V}$ is equivalent to the standard $L^2$-scalar product, yields $L^2$-stability.
  \end{proof}
\end{prop}

We have only shown $L^2$-stability. Unfortunately this is all we can hope for in
the relevant case of different transmission lines. In spite of this we will
sketch the following result:

\begin{prop}%
  \label{prop:5}
  For a graph $\mathscr{G} = (\mathscr{V},\mathscr{E})$ where all line parameters are equal: $(R_e,L_e,G_e,C_e) = (R,L,G,C)\ \forall e \in \mathscr{E}$ and an initial condition in $H^2$ (and sufficiently smooth coupling conditions) the above problem converges exponentially in the supremum norm.
\end{prop}

\autoref{prop:5} is proved by noting that we could choose $H^2$ for the domain of operator $\mathcal{A}$ in~\eqref{eq:78} without invalidating the proofs.
From the coupling conditions in~\eqref{IIBVP} we derive
\begin{equation*}
    \begin{aligned}
    \sum_{e \in \mathscr{E}_l} s(e,l)v_e^\prime (t,x_l) &= - R\sum_{e \in \mathscr{E}_l}s(e,l) i_e(t,x_l) - L \sum_{e \in \mathscr{E}_l}s(e,l)\dot{i}_e(t,x_l)\\
    &= -R i_l(t) - L \dot{i}_l(t) \eqqcolon v^\prime_l(t),
  \end{aligned}
\end{equation*}
and similarly
\begin{equation*}
    \begin{aligned}
    i^\prime_e &= i^\prime_f \ \forall\ e,f: x_v \in e,\\
    i^\prime_e &= i^\prime_g \ \forall g \in \mathcal{G}, e \in \mathscr{E}_g,
\end{aligned}
\end{equation*}
which means that the solution fulfills additional coupling conditions in this case.
Note that these calculations are only possible because the line parameters are equal and because the coupling conditions are assumed to be differentiable.
In this case, it is possible to mimic the stability proof above but for the enhanced Lyapunov function

\begin{equation*}
\mathcal{V}(u(t)) = \frac{1}{2} \sum_{e \in \mathscr{E}} \int_0^{\ell_e}C_e(u^1_e)^2 + L_e(u^2_e)^2 + \frac{1}{2} \sum_{e \in \mathscr{E}} \int_0^{\ell_e}C_e\left(\pd{u^1_e}{x}\right)^2 + L_e\left(\pd{u^2_e}{x}\right)^2 .
\end{equation*}

\section{From Telegrapher's equation to powerflow equations}
\label{sec:power-network}

Now that we have established that only the coupling conditions matter,
we examine periodic coupling conditions and their relation to the powerflow equations.
This means the functions $v_g$ and $i_l$ in \eqref{KCg} and \eqref{KCl} shall be of the form
\begin{equation*}
  \begin{aligned}
    v_g(t) &= \sum_{m \in \mathbb{Z}}\hat{V}_m^g e^{jm \omega t},\\
    i_l(t) &= \sum_{m \in \mathbb{Z}}\hat{I}_m^le^{jm \omega t} ,
  \end{aligned}
\end{equation*}
where quantities with a hat $\hat{\cdot}$ are prescribed at the node.
A smooth solution to~\eqref{IIBVP} with these coupling conditions is given by a Fourier series with coefficients of the form \eqref{eq:10},
where the coupling conditions translate to linear conditions on the complex voltage and current $V_m$ and $I_m$ for each transmission line.
The correct coefficients in this expression can be calculated by solving the linear system stemming from the coupling conditions.
As the system is linear, this can be done for each Fourier mode $m\omega$ separately.

The Fourier mode of the current at the ends of a line can be expressed as
\begin{equation*}
  \begin{pmatrix}
    I_m(0)\\
    I_m(\ell)
\end{pmatrix}
= -\frac{Y^0_m}{{\sinh(\gamma_m \ell)}}
\begin{pmatrix}
  - \cosh(\gamma_m \ell) & 1 \\
  -1 & \cosh(\gamma_m\ell) 
\end{pmatrix}
\begin{pmatrix}
  {V}^0_m \\
  {V}^\ell_m 
\end{pmatrix}.
\end{equation*}
We see that the ingoing end behaves differently from the outgoing end, because of an overall minus sign in the second component.
This is natural in the PDE setting as the sign of the current marks the direction of flow.
By using $I_m^0 = I_m(0)$ and $I_m^\ell = -I_m(\ell)$, we can transform the system into 
\begin{equation*}
    \begin{pmatrix}
    I_m^0\\
    I_m^\ell
  \end{pmatrix}
  = -\frac{Y^0_m}{{\sinh(\gamma_m \ell)}}
  \begin{pmatrix}
    - \cosh(\gamma_m \ell) & 1 \\
    1 & - \cosh(\gamma_m\ell) 
  \end{pmatrix}
  \begin{pmatrix}
    {V}^0_m \\
    {V}^\ell_m 
  \end{pmatrix}.
\end{equation*}
This makes the transmission line invariant under switching orientation
and positive currents mean currents leaving the node into the transmission line.
A further benefit is that we can omit the function $s$ from equation \eqref{eq:81}.
With this we come back to the whole network.
By the continuity condition on the voltage, i.e. \eqref{CC},
we do not need voltage Fourier modes for every transmission line,
but only for each vertex. Therefore we introduce the variables $V^r_m$
as the $m$th voltage Fourier mode at vertex $r$.

This is different for the current.
For a transmission between node $r$ and node $s$ we define the current leaving node $r$ in direction $s$ as
\begin{equation*}
  I_m^{rs} = -\frac{Y^{0,rs}_m}{{\sinh(\gamma_m^{rs} \ell^{rs})}}\left(- V_m^r\cosh(\gamma_m^{rs} \ell^{rs}) + V_m^s \right),
\end{equation*}
where $Y$, $\gamma$ and $\ell$ now depend on the transmission line.

The net current $I^r_m$ at node $r$ is then just the sum of all current leaving that node,
\begin{equation}
  \label{eq:36}
  I_m^r = \sum_{s \in \mathscr{E}} -\frac{Y^{0,rs}_m}{{\sinh(\gamma_m^{rs} \ell^{rs})}}\left(- V_m^r\cosh(\gamma_m^{rs} \ell^{rs}) + V_m^s \right)\ \forall r \in \mathscr{V},
\end{equation}
where we set $Y^{0,rs} = 0$ if node $r$ and $s$ are not connected.
These equations are then combined for the whole network, i.e., 
\begin{equation}
  \label{eq:70}
  I_m = Y_m V_m ,
\end{equation}
where $I_m,V_m \in \mathbb{R}^{\abs{\mathscr{V}}}$ are the vectors of net current and voltage at each node and $Y_m$ is the so-called admittance matrix.
One can set up admittance matrices for different power networks, where not every edge is a transmission line.
Yet in our case $Y_m$ is given as
\begin{equation*}
    \begin{aligned}
    (Y_m)_{ii}&=\sum_{j=1}^{\abs{\mathscr{V}}}\frac{Y^{0,ij}_m}{\tanh(\gamma_m^{ij}\ell^{ij})} ,\\
    (Y_m)_{ij}&=(Y_m)_{ji}=-\frac{Y^{0,ij}_m}{\sinh(\gamma_m^{ij}\ell^{ij})} .
\end{aligned}
\end{equation*}
An admittance matrix $Y_m$ is invertible, whenever not all of its row sums are zero (see~\cite[Theorem 1]{admittance_matrix_invertible}).
This is usually the case here, but we remark that it is possible to find a set of line parameters for which $Y_m$ is singular.

That means the coupling conditions in the time periodic setting are usually solvable if all nodes are loads.
If there are generator nodes, \eqref{eq:70} can be reduced by inserting the values of $\hat{V}_m^g$ and
shifting them to the other side by altering $I_m$.
Concretely, we define the vector $\hat{V}_m$ by
\begin{equation*}
    \hat{V}_m^i =
  \begin{cases}
    V_m^i, \quad \text{ if the $i$th node is a generator}\\
    0, \quad \quad \text{ else}
  \end{cases}
\end{equation*}
and
\begin{equation*}
    \hat{I}_m = I_m-Y_m\hat{V}_m .
\end{equation*}
Then, we strike all lines (and columns respectively) from $\hat{I}_m$ and $Y_m$ that correspond to generator nodes and arrive at a new equation
\begin{equation*}
  \tilde{I}_m = \tilde{Y}_m \tilde{V}_m,
\end{equation*}
in which $\tilde{Y}_m$ is usually\footnote{Again by carefully choosing line parameters it is possible to make the resulting $\tilde{Y}_m$ singular.} still invertible and only load nodes are left.
Therefore apart from pathological counter-examples the coupling conditions,
\begin{equation}
  \label{eq:69}
  \begin{aligned}
    \hat{V}_m^g&=V^g_m \ \forall \ g \in \mathcal{G},\\
    \hat{I}_m^l &= \sum_{s \in \mathscr{E}} -\frac{Y^{0,ls}_m}{{\sinh(\gamma_m^{ls} \ell)}}\left(- V_m^l\cosh(\gamma_m^{ls} \ell) + V_m^s \right)\ \forall l \in \mathcal{L},
  \end{aligned}
\end{equation}
which are again linear in the complex voltages, are solvable by functions defined by equations \eqref{eq:2} through \eqref{eq:10}.

\subsection{Power coupling conditions}%
\label{sec:power-coupl-cond}
In power networks we usually want to satisfy a power demand instead of a current demand.

Electric power at a node $v$ is given by the product of voltage and the net current,
\begin{equation*}
    \begin{aligned}
    p_v(t) &= v_v(t)\cdot i_v(t) 
    &= \sum_{m,k\in \mathbb{Z}} V_m I_k e^{j (m+k)\omega t} 
    &= \sum_{m,k \in \mathbb{Z}} V_{m-k} I_k e^{j m\omega t}
    &\eqqcolon \sum_{m \in \mathbb{Z}} P_m e^{j m\omega t} . 
  \end{aligned} 
\end{equation*}
One way to replace the linear coupling conditions~\eqref{eq:69} would be to prescribe $\hat{P}_m^l $ at every load node.
In the case where $v_g$ and $p_l$ are Fourier polynomials, that is there is $M \in \mathbb{N}_0$, such that $\hat{P}_m^l,V_m^g = 0$ for $\abs{m}> M$,
this yields $\abs{\mathscr{V}}(2M+1)$ equations for just as many unknown voltage Fourier modes%
\footnote{Up to now these are complex equations for complex Fourier modes but assuming $V_m = V_{-m}^*$ yields $\abs{\mathscr{V}}(2M+1)$ real equations for as many real unknowns.}.
So depending on the concrete values of $\hat{P},\hat{V}$ the system might be solvable.
Unfortunately, except for small examples it is hard to decide a priori whether a given set of power coupling conditions
is feasible, see for example \cite{infeasible_powerflow}.

The mixing of different Fourier modes can make the resulting coupling conditions quite cumbersome.
Luckily in applications a strong simplification is done:
It is assumed that all inputs are sinusoidal,
meaning that the only non-vanishing Fourier modes are those for $m \in \set{\pm 1}$.

We use $v,i \in \mathbb{R}$ to express $V_{-1} = V_1^*$ and define  $V \coloneqq \abs{V}e^{i\phi}$, similarly $I \coloneqq \abs{I}e^{i\psi},$ and find
for the power
\begin{equation*}
    \begin{aligned}
    p(t) &= \left(Ve^{i\omega t}+V^*e^{-i\omega t}\right)\left(Ie^{i\omega t}+I^*e^{-i\omega t}\right)\\
    &=VI^* +V^*I + VIe^{2i\omega t} +V^*I^*e^{-2i\omega t}\\
    &=\Re(VI^*) + \Re(VIe^{2i\omega t}) \\
    &=\Re(VI^*)\left(1 + \cos(2(\omega t+\phi))  \right) + \Im(VI^*) \sin(2(\omega t +\phi)) .
  \end{aligned}
\end{equation*}
This suggests to define the so-called \emph{complex apparent power} $S = VI^*$, its real part,
the \emph{real power} $P = \Re(S)$ and its imaginary part, the reactive power $Q = \Im(S)$,
so that at the node $v$ it holds:
\begin{equation*}
    p_v(t) = P_v\left(1 + \cos(2(\omega t+\phi))  \right) + Q_v \sin(2(\omega t +\phi)) .
\end{equation*}
To define the powerflow equations we use the admittance matrix $Y$ from \eqref{eq:36}:
\begin{equation*}
    I^* = (Y V)^*,
\end{equation*}
of which we write the real and imaginary part as
\begin{equation*}
    Y = G +iB.
  \end{equation*}
  For the complex apparent power this means
  \begin{equation*}
  S = \diag(V)Y^*V^* .
\end{equation*}
For every node $s$ we then find the nonlinear system of equations
\begin{equation*}
    \begin{aligned}
    P_{s}  & = \sum_{i=1}^N \abs{V_s}\abs{V_i}(G_{si}\cos(\phi_s-\phi_i)+B_{si}\sin(\phi_s-\phi_i)),\\
    Q_{s}  & = \sum_{i=1}^N   \abs{V_s}\abs{V_i}(G_{si}\sin(\phi_s-\phi_i)-B_{si}\cos(\phi_s-\phi_i)),
  \end{aligned}
\end{equation*}
which are the so-called \emph{powerflow equations}.

The coupling conditions for our network problem are then usually chosen as follows:
\begin{itemize}
\item All load nodes prescribe $P$ and $Q$.
\item All generators (with one exception) prescribe $\abs{V}$ and $P$.
\item The last generator, known as the \emph{slack bus}, prescribes $\abs{V}$ and $\phi$.
\end{itemize}

These conditions lead to a system of nonlinear equations with two equations for each node
to determine the two nodal variables $\abs{V}$ and $\phi$. For more information, we refer to~\cite{Frank16,Andersson15book}. 

Unfortunately it is unclear how the system of powerflow equations can be transformed back into the PDE setting as it has an inherently periodic time dependence.
Therefore, in order to examine the stability properties of the system in the following section,
we will treat nodal voltages and currents only as linear coupling conditions in~\eqref{IIBVP}.
We prescribe the net current at load nodes and the voltage at generators.
The justification for this is, that ``voltage is pushed'' (set by the voltage source) and ``current is pulled'' (drawn by the load).
\section{Numercial scheme}%
\label{sec:numercial-scheme}
We aim to mimic the behavior of the analytical solution of~\eqref{IIBVP} by a numerical approximation.
We will first treat a single line, define the numerical scheme there and then
prescribe numerical boundary conditions.

We approximate a solution of~\eqref{eq:6} by an approximation and compute the discrete solution
via a splitting scheme. The numerical approximation is based on equation~\eqref{eq:14}, where
\begin{equation*}
  \Lambda =\diag( \lambda_1, \dots , \lambda_N, -\lambda_1, \dots ,
  -\lambda_N)
\end{equation*}
is diagonal.
The balance law is split into an ODE part
\begin{equation}
  \label{eq:31}
  \xi_t = -B\xi
\end{equation}
and a PDE part
\begin{equation}
  \label{eq:32}
  \xi_t +\Lambda \xi_x=0.
\end{equation}
Therefore we choose the discretization $t_n = n\Delta t$, $n = 0\dots N$, $x_j = j \Delta x$, $j = 0\dots j$ and
\begin{equation*}
    {(\xi^{\pm})}^n_{j} = \xi^{\pm}(t_n,x_j) .
\end{equation*}
Further, we define the lattice constant as
\begin{equation*}
    \Gamma = \frac{\Delta t}{\Delta x} .
\end{equation*}
As splitting we use Strang splitting (see~\cite[17.4 Strang splitting]{leveque2002finite}), meaning we make a half-step of equation~\eqref{eq:31}, a full step of~\eqref{eq:32} and another half-step of~\eqref{eq:31}.

Equation~\eqref{eq:31} is a linear ODE with constant coefficients, so we can solve it exactly up to machine
precision with only constant computational costs by setting the ODE scheme to
\begin{equation*}
    \begin{aligned}
    \xi^*_j &= {(\ODE(\Delta t)\cdot \xi)}_j =  M\xi_j^n ,\\
    M &= \exp(-B\Delta t).
  \end{aligned}
\end{equation*}
As Strang splitting is of second order, we cannot hope to get higher accuracy and therefore choose a scheme for the conservation law~\eqref{eq:32} which is at most second order,
namely a flux-limited Lax-Wendroff scheme (explained for example in~\cite[6.12 TVD limiters]{leveque2002finite}) that for a single component is given for right and left moving waves as: 
\begin{equation*}
    \begin{aligned}
    \hat{\xi}^+_j &= \xi^+_j - \Gamma \abs{a} (\xi^+_{j}-\xi^+_{j-1}) \\
    &- \frac{1}{2}\abs{a}\Gamma(1-\abs{a}\Gamma)\left(\phi(\theta_{j+\frac{1}{2}}^+)(\xi^+_{j+1}-\xi^+_{j})-\phi(\theta_{j-\frac{1}{2}}^+)(\xi^+_{j}-\xi^+_{j-1})\right),\\
    \hat{\xi}^-_j &= \xi^-_j + \Gamma \abs{a} (\xi^-_{j+1}-\xi^-_{j}) \\
    &- \frac{1}{2}\abs{a}\Gamma(1-\abs{a}\Gamma)\left(\phi(\theta_{j+\frac{1}{2}}^-)(\xi^-_{j+1}-\xi^-_{j})-\phi(\theta_{j-\frac{1}{2}}^-)(\xi^-_{j}-\xi^-_{j-1})\right),
  \end{aligned}
\end{equation*} 
where $a$ is the wave speed and we define as usual
\begin{equation*}
    \begin{aligned}
    \theta_{j+\frac{1}{2}}^+ &= \frac{\xi^+_{j\phantom{{}+1}}-\xi^+_{j-1}}{\xi^+_{j+1}-\xi^+_{j\phantom{{}+1}}},\\
    \theta_{j+\frac{1}{2}}^- &= \frac{\xi^-_{j+2}-\xi^-_{j+1}}{\xi^-_{j+1}-\xi^-_{j\phantom{{}+1}}}.
  \end{aligned}
\end{equation*}
As limiter we choose the minmod-limiter,
\begin{equation*}
    \phi(x) =
  \begin{cases}
    0 \text{ for }\phantom{0<{}}x<0,\\
    x \text{ for }0<x<1,\\
    1 \text{ for }1<x,
  \end{cases}
\end{equation*}
which will become important for the numerical Lyapunov stability.
As with any explicit finite volume scheme a necessary condition for stability is the CFL condition
\begin{equation}
  \label{eq:8}
  \Gamma \abs{a} \leq 1,
\end{equation}
which relates the lattice constant and the wave speed.  The quantity $\Gamma \abs{a}$ is also called the \emph{CFL number}.

We write the scheme (although it is not linear) as
\begin{equation*}
    \hat{\xi}_j = {(\PDE(\Delta t)\xi)}_j,
\end{equation*}
so that the full split scheme can be written as
\begin{equation*}
    \xi^{n+1} = \ODE\left(\frac{\Delta t}{2}\right)\circ \PDE(\Delta t)\circ \ODE\left(\frac{\Delta t}{2}\right)\cdot\xi^n .
\end{equation*}
Having defined the scheme we come to the two main properties of it. We will now prove that the scheme is 
\begin{itemize}
\item total-variation-diminising and
\item Lyapunov stable.
\end{itemize}
As numerical Lyapunov function we choose 
\begin{equation}
  \label{eq:37}
  \mathcal{V}(\xi) = \frac{1}{2} \sum_j\left( {(\xi_j^+)}^2+ {(\xi_j^-)}^2\right) = \frac{1}{2}\sum_j {(\xi_j)}^T\xi_j,
\end{equation}
which is just the squared discrete $L^2$-norm of the solution and hence corresponds directly to the analytical Lyapunov function.

As definition of total variation in two dimensions we choose the sum of the componentwise total variations:
\begin{equation}%
  \label{eq:3}
    \TV(\xi) = \sum_j \left(\,\abs{\xi^+_{j+1}-\xi^+_{j}}+\abs{\xi^-_{j+1}-\xi^-_{j}}\right) = \sum_j \norm{\xi_{j+1}-\xi_j}_1 .
\end{equation}

Now we are able to discuss the following results:
\begin{prop}%
  \label{prop:2}
  The ODE-scheme~\eqref{eq:31} is Lyapunov stable and total-variation-diminishing.
\end{prop}
For the proof we need a short lemma.
\begin{lemma}%
  \label{lemma:1}
  There holds
  \begin{enumerate}
  \item\label{item:3}The matrix $M= \exp(-B\Delta t)$  has the structure
    \begin{equation}
      \label{eq:39}
      M =
      \begin{pmatrix}
        r&s\\s&r
      \end{pmatrix}.
    \end{equation}
  \item\label{item:4}Its eigenvalues are given by $m_1 = r+s$, $m_2=r-s$ and fulfill $0<m_1,m_2<1$.
  \item\label{item:5}There holds $r>0$.
  \end{enumerate}
  \begin{proof}
    Consider the base change matrix
    \begin{equation*}
            P =
      \begin{pmatrix}
        0&1\\1&0
      \end{pmatrix} .
    \end{equation*}
    Note that a matrix $A$ fulfills
    \begin{equation*}
            P^{-1}AP = A,
    \end{equation*}
    if and only if $A$ has the form~\eqref{eq:39}.
    Therefore there holds $P^{-1}BP = B$. But $\exp$ commutes with change of basis and therefore
    \begin{equation*}
            P^{-1}MP = \exp(-P^{-1}BP \Delta t) = \exp(-B\Delta t) = M,
    \end{equation*}
    proving~\ref{item:3}.
    The eigenvalues of $M$ can be computed from~\eqref{eq:39} and are $m_{1,2} = r\pm s$.
    So they are also given by $m_{1,2} = \exp(-b_{1,2}\Delta t)$ where $b_{1,2}$ are the eigenvalues of $B$.
    Assertation \ref{item:4} follows because $b_{1,2}>0$.
    Lastly $r>= \frac{m_1+m_2}{2}>0$.
  \end{proof}
\end{lemma}
Next we come to the
\begin{proof}[Proof of \autoref{prop:2}]
  For the total variation after an ODE step we find
  \begin{equation*}
        \begin{aligned}
      \TV(\ODE(\Delta t)\xi) &=  \sum_j \norm{M\xi_{j+1}-M\xi_j}_1\\
      &\leq \norm{M}_1 \sum_j \norm{\xi_{j+1}-\xi_j}_1 \\
      &= \left(\,\abs{r}+\abs{s}\right) \sum_j \norm{\xi_{j+1}-\xi_j}_1\\
      &\leq \max(r+s,r-s)\sum_j \norm{\xi_{j+1}-\xi_j}_1\\
      &\leq \sum_j \norm{\xi_{j+1}-\xi_j}_1 = \TV(\xi),
    \end{aligned}
  \end{equation*}
  due to \autoref{lemma:1}, specifically, the second inequality is due to assertion~\ref{item:5} and the last one due to assertion~\ref{item:4}.

  For the Lyapunov function after an ODE step we get
  \begin{equation*}
        \begin{aligned}
      \mathcal{V}(\ODE \xi) &= \frac{1}{2} \sum_j {(M\xi_j)}^{T}M\xi_j\\
      &=\frac{1}{2} \sum_j {(\xi_j)}^{T}M^{T}M\xi_j\\
      &\leq \frac{1}{2} \max(m_1^2,m_2^2)\sum_j {(\xi_j)}^{T}\xi_j\\
      &\leq \frac{1}{2} \sum_j {(\xi_j)}^{T}\xi_j = \mathcal{V}(\xi),
    \end{aligned}
  \end{equation*}
  again according to \autoref{lemma:1}.
\end{proof}

On the PDE side we find similar results.
\begin{prop}%
  \label{prop:3}
  The PDE-scheme is Lyapunov stable and total-variation-diminishing.
  \begin{proof}
    For the claim of TVD we note that the total variation~\eqref{eq:3} is given as the sum of the variations in both components.
    As the PDE-scheme does not mix the components, it is overall total variation diminishing, if it is TVD each of its components. 
    The TVD property of the minmod-limiter for the scalar case is a standard result found e.g.\ in~\cite[6.12 TVD Limiters]{leveque2002finite}.

    The Lyapunov function~\eqref{eq:37} is also given as a sum of the components, namely the sum of the $L^2$-norms squared.
    Therefore we examine $L^2$-stability of the PDE-scheme for the linear advection equation (which is one component of~\eqref{eq:32}), given by
    \begin{equation}
      \label{eq:47}
      u_t +au_x=0.
    \end{equation} with the wave speed $a$. 
    The pair
    \begin{equation*}
            E(u) = \frac{1}{2}u^2 , \ Q(u) = \frac{a}{2}u^2
    \end{equation*}
    is an entropy-entropy-flux pair for~\eqref{eq:47}.
    According to~\cite[Theorem 4.1]{tadmore1987} there is an entropy-conserving three-point scheme for this entropy.
    This scheme turns out to be the so-called naive scheme with numerical flux
    \begin{equation*}
            F_{j+\frac{1}{2}} = \frac{a}{2}(u_j+u_{j+1})
    \end{equation*}
    and numerical viscosity
    \begin{equation*}
            Q^{\text{naive}}_{j+\frac{1}{2}} = 0.
    \end{equation*}
    Now according to~\cite[Theorem 5.2]{tadmore1987} a scheme with viscosity $Q_{j+\frac{1}{2}}$ is entropy-stable if
    \begin{equation*}
            Q_{j+\frac{1}{2}}\geq Q^{\text{naive}}_{j+\frac{1}{2}} = 0
    \end{equation*}
    for all $j$.
    The viscosity of the minmod scheme is given by
    \begin{equation*}
            Q^{mm}_{j+\frac{1}{2}} = \abs{a}\Gamma\left(1-\left[1-\abs{a}\Gamma)\phi(\theta^{\sign(a)}_{j+\frac{1}{2}})\right]\right) \geq 0,
    \end{equation*}
    because of the CFL condition \eqref{eq:8} and $\phi\leq 1$.
  \end{proof}
\end{prop}

\subsection{Coupling conditions}
\label{sec:coupling-conditions}
Now that we have established the stability of the scheme in the interior of the
transmission line, we need to define the coupling conditions.
Therefore we consider a node $n$ with a set $\mathscr{E}_n$ of transmission lines, that
all start at it and reformulate the combined Telegrapher's equations on all
these lines.
We define the vector and the matrices
\begin{equation*}
  \begin{aligned}
        \xi &=
    \begin{pmatrix}
      \xi_1^+& \dots & \xi_N^+& \xi_1^-& \dots & \xi_N^-
    \end{pmatrix}^T,\\
    \Lambda &=\diag( \lambda_1, \dots , \lambda_N, -\lambda_1, \dots ,
    -\lambda_N)
    ,\\
    B &=
    \begin{pmatrix}
      \alpha & \beta\\
      \beta & \alpha
    \end{pmatrix},\\
    \alpha &= \diag(a_1, \dots , a_N),\\
    \beta &= \diag(b_1, \dots , b_N)
  \end{aligned}
\end{equation*}
and also
\begin{equation*}
    \xi^\pm =
  \begin{pmatrix}
    \xi_1^\pm&
    \dots &
    \xi_N^\pm&
  \end{pmatrix}^T, \
  \Lambda^\pm =\pm\diag( \lambda_1,  \dots ,  \lambda_N) ,
\end{equation*}
which leads to the PDE
\begin{equation*}
    \xi_t +\Lambda \xi_x + B\xi=
  \begin{pmatrix}
    \xi^+\\
    \xi^-
  \end{pmatrix}_t
  +
  \begin{pmatrix}
    \Lambda^+&\\
    &\Lambda^-
  \end{pmatrix}
  \begin{pmatrix}
    \xi^+\\
    \xi^-
  \end{pmatrix}_x +
  \begin{pmatrix}
    \alpha & \beta\\
    \beta & \alpha
  \end{pmatrix}
  \begin{pmatrix}
    \xi^+\\
    \xi^-
  \end{pmatrix}
  = 0 .
\end{equation*}

For line $e \in \mathscr{E}_n$ we define $c_e = \sqrt{\frac{L_e}{C_e}}$. If the node is a load with prescribed
current $i_n$, we have the following conditions:
\begin{equation*}
    \begin{aligned}
    \sum_{i = 1}^Ni_e  = i_n  &\Leftrightarrow    \sum_{e = 1}^N\xi_e^+  = -\sum_{e = 1}^N  \xi_e^- + i_n,\\
    v_e = v_f\ \forall \ 1\leq e,f \leq N    &\Leftrightarrow    c_e(\xi_e^+ - \xi_e^-) = c_f(\xi_f^+ - \xi_f^-).
  \end{aligned}
\end{equation*}
To reformulate them we define the matrices
\begin{equation*}
    M =
  \begin{pmatrix}
    1 & &\dots  & 1\\
    c_1& -c_2 \\
    &\ddots&\ddots\\
    &&c_{N-1}&-c_N
  \end{pmatrix}, \
  S = \begin{pmatrix}
    -1\\
    &1\\
    &&\ddots\\
    &&&1
  \end{pmatrix}, \
  \tilde{i}_n = \begin{pmatrix}
    i_n \\
    0\\
    \vdots \\
    0
  \end{pmatrix},
\end{equation*}
which yields the conditions\footnote{Note that $M$ is invertible, as its determinant is given by
  $ = -{(-1)}^N (\prod_{i=1,}^N c_i) (\sum_{i=1}^N\tfrac{1}{c_i})$ and $c_i>0$.} for the ghost cell values $\xi_0$:
\begin{equation}
  \label{eq:65}
  M\xi_0^+ =  SM\xi_0^- + \tilde{i}_n \Leftrightarrow \xi_0^+ = M^{-1}SM\xi_0^-
  + M^{-1}\tilde{i}_n \Leftrightarrow\ \xi_0^+ = U\xi_0^- + M^{-1}\tilde{i}_n .
\end{equation}
Note that $\norm{U}_2 = 1$, so for $i_n = 0$ no energy is injected at the node.
We still must choose $\xi^-_0$ and will do so by linear extrapolation.
In principle one should use stable extrapolation (see~\cite{nonsmooth_extrapolation}) to guard against non-smooth solutions
but this seems to be unnecessary with high enough $R$ and $G$.

For the mimmod scheme we need a second ghost cell $\xi_{-1}$.
In principle it should be possible to use the time derivative of the coupling condition as in~\cite{Banda_2016} or~\cite{Borsche_2014}, but in our setting it is much easier to simply use a first order extrapolation through the first couple of inner cell values and the newly-computed boundary value from~\eqref{eq:65}.

Contrary, the generator coupling conditions are easier to determine as they essentially decouple
different lines.
When we prescribe the voltage at a node, this acts as a regular boundary
condition for each attached line and can be treated as such.
Alternatively we can write the generator coupling conditions similar to equation~\eqref{eq:65},
\begin{equation}
\label{eq:5}  
    \xi^+_0 = \xi^-_0+  {\diag(c_1,\dots,c_N)}^{-1} {(1,\dots,1)}^{T} v_n.
\end{equation}
Regarding the stability of coupling conditions~\eqref{eq:65} and~\eqref{eq:5} we make the following remark.

\begin{remark}
We note that any numerical entropy production at the nodes is countered by the decay in the ODE part, whenever $R$ and $G$ are great enough.
Unfortunately we suppose that these coupling conditions are unstable if $R$ and $G$ are small enough. A rigorous analysis will be subject of future work.
\end{remark}

\section{Numerical examples}\label{sec:numerical-examples}

For the computations we wrote a package\footnote{The package is free software and be used and modified by anyone. It can be found under \url{https://bitbucket.org/efokken/telegraph_numerics}} for the Julia programming language\footnote{see \url{https://julialang.org/}} (see~\cite{julialang}).

We will treat two examples, in which we compare exact periodic solutions, that correspond to powerflow solutions and numerical solutions computed with the splitting scheme from section \ref{sec:numercial-scheme}. In addition we consider the behavior of the numerical Lyapunov function.

As we compare the splitting scheme to the powerflow solution, the coupling conditions in the node are given by sinusoidal functions, namely by
\begin{equation*}
    \begin{aligned}
    v(t) &= \Re(V e^{\omega j t}),\\
    i(t) &= \Re(I e^{\omega j t}).\\
  \end{aligned}
\end{equation*}
Therefore we provide only the quantity $I$ and $V$ for each node and of course $\omega$ for the whole network.
The first example under consideration is a single transmission line between two nodes, one of which supplies voltage, while the other draws current.

\begin{figure}[ht]
  \centering
  \begin{subfigure}[t]{0.49\linewidth}
    \centering
    \includegraphics[width=\textwidth]{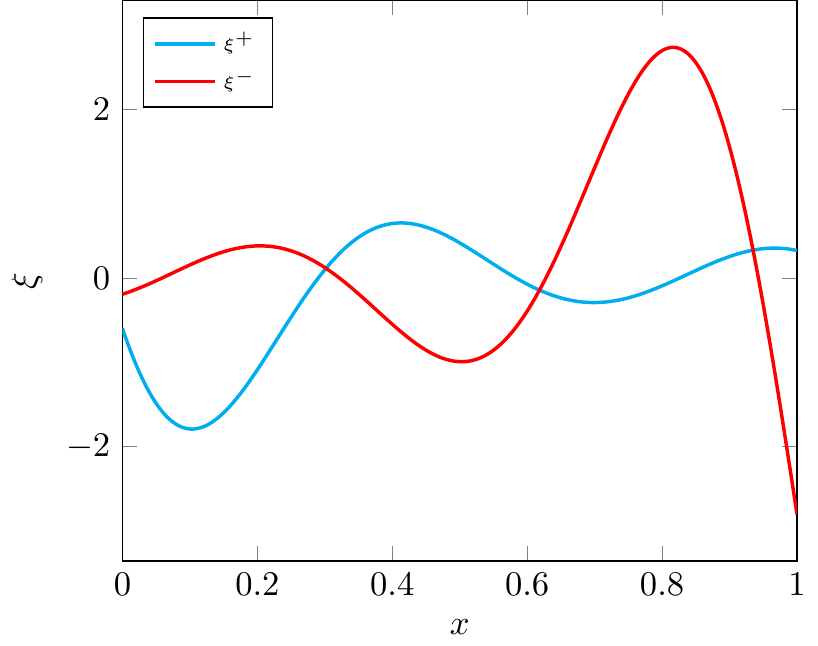}
    \caption{Solution at $t=1.0$.}
    \label{fig:sol1}
  \end{subfigure}
  \begin{subfigure}[t]{0.49\linewidth}
    \centering
    \includegraphics[width=\textwidth]{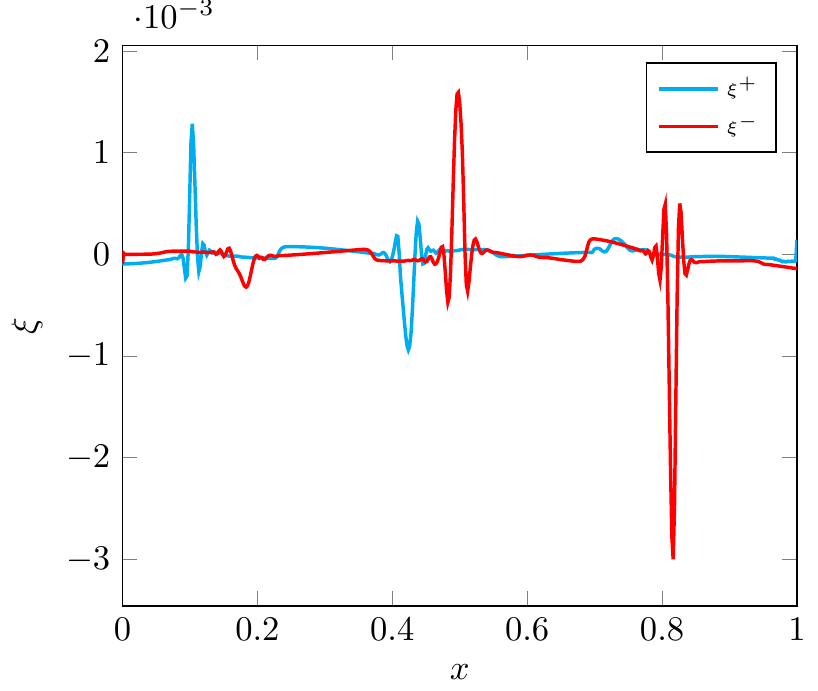}
    \caption{Error of the solution at $t=1.0$.}
    \label{fig:error}
  \end{subfigure}
  \caption{The single line example.}
\end{figure}

The parameters of this network are given in \autoref{tab:example1} and the parameters of the periodic solution in \autoref{tab:example1periodic}.
\begin{table}[ht]
  \centering
  \caption{Parameters of the line, in addition $\omega = 4.0$.}
  \begin{subtable}{.49\textwidth}
    \centering
    \caption{Line parameters.}
    \begin{tabular}{lr}
      \toprule
      $R$& $4.0$\\
      $L$& $6.0$\\
      $G$& $2.0$\\
      $C$& $1.0$\\
      $\ell$& $1.0$\\
      \bottomrule
    \end{tabular}\label{tab:example1}
  \end{subtable}
  \begin{subtable}{.49\textwidth}
    \centering
    \caption{Node parameters.}
    \begin{tabular}{lr}
      \toprule
      $V_\text{start}$ & $5.0+3.0j$\\
      $I_\text{end}$ & $2.0+5.0j$\\
      \bottomrule
    \end{tabular}
    \label{tab:example1periodic}
  \end{subtable}
\end{table}
We integrate this problem from $t=0$ to $t=1.0$ with space discretization $\Delta x =2^{-9}$.
As CFL number we choose $0.8$.
The numerical solution at $t=1.0$ can be found in \autoref{fig:sol1}.
To compute the error we first provide the analytical solution expressed in $\xi^{\pm}$.
The analytical solution on every line is given by
\begin{equation*}
    \begin{aligned}
    v(x,t) = \Re(V(x)e^{\omega j t}),\\
    i(x,t) = \Re(I(x)e^{\omega j t}),\\
  \end{aligned}
\end{equation*}
where $V,I$ are taken from equation~\eqref{eq:10} (with $m=\pm 1$).
As voltages at the ends of the line we take the solution to the linear system~\eqref{eq:69}.
With these we compute the analytical expressions for $\xi^{\pm}$ as
\begin{equation*}
    \xi^{\pm}(x,t) = \frac{1}{2}\left(i(x,t) \pm \sqrt{\frac{C}{L}}v(x,t) \right)
\end{equation*}
on every line. 
For an analytical solution $\zeta$ and a numerical solution $\xi$ the maximal error of the numerical solution on a network over a time $T = J\Delta t$ is then given by
\begin{equation}
  \label{eq:89}
  \Delta \xi=\max_{j=1:J}\max_{e \in \mathscr{E}}(\max_{k=1:N_e}(\abs{\xi^{+}_{kj}-\zeta^{+}_{kj}},\abs{\xi^{-}_{kj}-\zeta^{-}_{kj}})),
\end{equation}
where $\mathscr{E}$ is again the set of edges and $N_e$ is the number of spatial steps on edge $e$ and $J$ is the total number of time steps.
Here, $\zeta_{kj} = \zeta(k\Delta x,j\Delta t)$ is just the analytical solution evaluated at cell centers.

The time-dependent error over the whole line at time $t=1$ is shown in \autoref{fig:error} for illustrating its usual appearance.
The spikes in the error chart come from the minmod-scheme and lie exactly at the extrema of the solution.
We remark that the boundary conditions are not introduce no further error.

\begin{figure}[ht]
  \centering
  \begin{tikzpicture}
    \node[draw, circle, fill=green!50] (N2) at (3.0,-1.0) {N2};
    \node[draw, circle, fill=blue!50] (N1) at (4.0,0.0) {N1};
    \node[draw, circle, fill=green!50] (N3) at (5.0,-1.0) {N3};
    \node[draw, circle, fill=green!50] (N4) at (4.0,1.2) {N4};

    \draw[thick] (N2) edge (N1);
    \draw[thick] (N3) edge (N1);
    \draw[thick] (N4) edge (N1);
  \end{tikzpicture}
  \caption{A network with three generators and one load.}
  \label{fig:network}
  \end{figure}
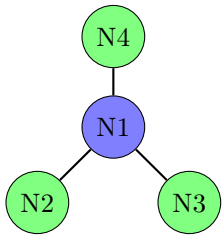

A second example is given by the network shown in \autoref{fig:network}, where the node in the middle is a current supplying node where all transmission lines start.
The line parameters and coupling conditions are given in table \autoref{tab:network}.  The network is also integrated with a spatial stepsize of $\Delta x =2^{-9}$ and again a CFL number of $0.8$.
For this network we also examine the order of convergence.
Therefore, we solve the problem with stepsizes $\Delta x_i = 2^{-i}$ for $1\leq i \leq 10$ and compute the corresponding maximal errors $\Delta \xi_i$ (see~\eqref{eq:89}) as well as the $i$th estimate of the convergence order by
\begin{equation*}
    c_i = \frac{\frac{\Delta \xi_i}{\Delta \xi_{i-1}}}{\frac{\Delta x_i}{\Delta x_{i-1}}}.
\end{equation*}

In \autoref{tab:conv-order} we see that the order of convergence approaches $2$ with the pure Lax-Wendroff scheme and stays below $2$ for the minmod-scheme as can be expected due to the order reduction at extrema.

\begin{table}[ht]
  \centering
  \caption{Parameters of the network, in addition $\omega = 4.0$.}
  \begin{subtable}{.49\textwidth}
    \centering
    \caption{Line parameters.}
    \begin{tabular}{rrrr}
      \toprule
      line & $1\rightarrow 2$&$1\rightarrow 3$&$1\rightarrow 4$\\
      \midrule
      $R$& 2.0&3.0&1.0\\
      $L$& 6.0&6.0&9.0\\
      $G$& 2.0&1.0&2.0\\
      $C$& 1.0&1.0&1.0\\
      $\ell$& 2.0&2.0&2.0\\
      \bottomrule
    \end{tabular}
  \end{subtable}
  \begin{subtable}{.49\textwidth}
    \centering
    \caption{Node parameters.}
    \begin{tabular}{rrr}
      \toprule
      node&type&coupling condition\\
      \midrule
      N1&current&10.0+3.0j\\
      N2&voltage&4.0+4.0j\\
      N3&voltage&2.0+5.0j\\
      N4&voltage&3.0+6.0j\\
      \bottomrule
    \end{tabular}
  \end{subtable}
  \label{tab:network}
\end{table}
	
\begin{figure}[ht]
  \centering
  \begin{subfigure}[t]{0.49\linewidth}
    \includegraphics[width=\textwidth]{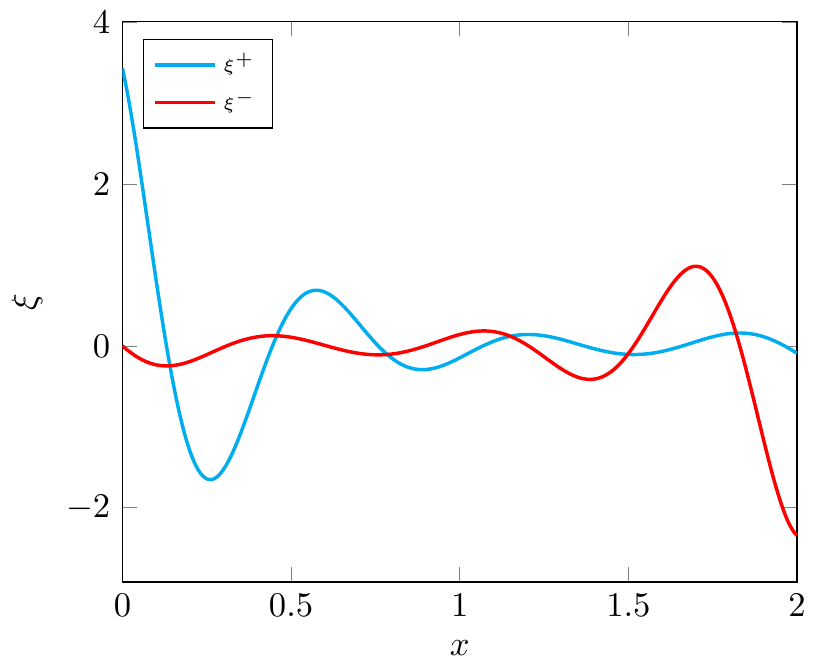}
    \caption{Solution on the line from N1 to N2 at $t = 3$.}%
    \label{fig:line12}
  \end{subfigure}
  \begin{subfigure}[t]{0.49\linewidth}
    \includegraphics[width=\textwidth]{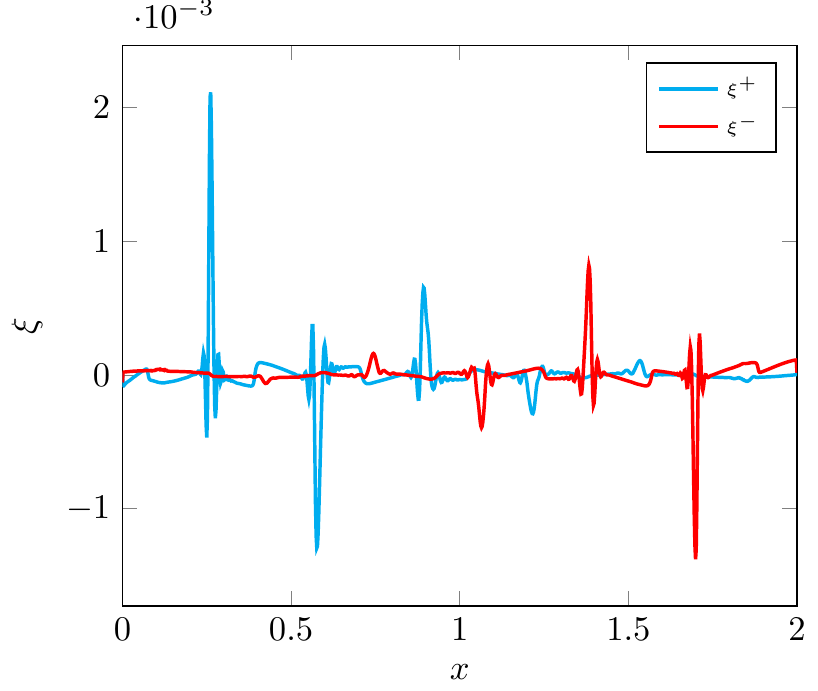}
    \caption{Error on the line from N1 to N2 at $t = 3$.}%
    \label{fig:error12}
  \end{subfigure}
  \caption{The network example.}
\end{figure}

\begin{table}[ht]
  \centering
  \caption{Order estimates of split schemes with different PDE schemes.}
  \begin{tabular}{rS[round-mode=places,round-precision=3]S[round-mode=places,round-precision=3]}
    \toprule
    {$i$} & {Lax Wendroff} & {Minmod}\\
    \midrule
    $2$ & 0.5743230046610167& 0.5837889926765258\\
    $3$ & 1.1369906850649314& 0.8464848165843416\\
    $4$ & 1.7542954324827598& 1.303766124009856 \\
    $5$ & 1.986301121408866 & 1.2359500844726483\\
    $6$ & 1.9137043365018416& 1.14530589269782  \\
    $7$ & 1.9545064711110431& 1.212455606493242 \\
    $8$ & 1.976492993708259 & 1.245011121236105 \\
    $9$ & 1.9880835112291628& 1.2930895877839055\\
    $10$& 1.993974711479301 & 1.3693499161473541\\
    \bottomrule
  \end{tabular}%
  \label{tab:conv-order}
\end{table}
		  Plots of the line solutions and errors are similar to those of the single line and can be found in \autoref{fig:line12} and \autoref{fig:error12}.  Different parameters do not change the pictures much, which is why we do not show further images.
  
Lastly we compare the numerical~\eqref{eq:37} and analytical Lyapunov functions .
The analytical Lyapunov function is given by
\begin{equation*}
    \mathcal{V}(t) = \mathcal{V}(0) \exp\left\{-2\min_{e \in \mathscr{E}} \left(\min\left( \frac{R_e}{L_e}, \frac{G_e}{C_e}\right)\right) t\right\}.
\end{equation*}
For the comparison we use the same network and the same initial data, but choose $\tilde{i}_l = 0$ and $v_g = 0$ for load and generator nodes respectively. 
With the numerical scheme we compute a solution and evaluate its Lyapunov function.
The analytical Lyapunov function in this case is given by $\mathcal{V}(0) \exp(-\frac{2}{3}t)$.

\begin{figure}[ht]
  \centering
  \includegraphics[width=\textwidth]{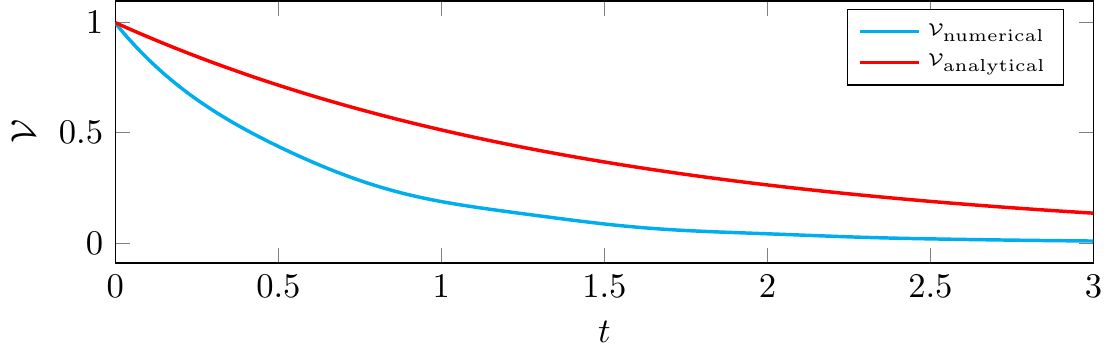}
  \caption{Comparison of the analytical estimate and the computed numerical Lyapunov function, both normed to $1$.}%
  \label{fig:lyapunov_net}
\end{figure}

In \autoref{fig:lyapunov_net} we observe the estimate and the actually computed Lyapunov function.

Lastly we want to illustrate the severe restriction, the CFL number places on utilizing an explicit scheme for realistic power networks.
Therefore we take some sample line parameters from~\cite[Table 3 on page 26]{freileitung_kabel}, namely parameters for an aerial transmission line for a voltage of $\SI{380}{\kilo \volt}$:
\begin{equation}
  \label{eq:9}
  \begin{aligned}
    R&= \SI{0.028}{\milli \ohm \per \kilo \meter},\\
    L&= \SI{0.8}{\milli \henry \per \kilo \meter},\\
    G&= \SI{15}{\nano \siemens \per \kilo \meter},\\
    C&= \SI{14}{\nano \farad \per \kilo \meter}
  \end{aligned}
\end{equation}
and compute the wave speed given in~\eqref{eq:1} from them,
\begin{equation*}
  \lambda = \sqrt{LC}^{-1} \approx \SI{298807}{\kilo \meter \per \second} \approx \SI{3e8}{\meter \per \second} ,
\end{equation*}
which is nearly the speed of light.  The wave speed $\lambda$ for a buried cable still fulfills $\lambda >\SI{1e8}{\meter \per \second}$.
Therefore for a rather coarse spatial step size of $\Delta x = \SI{1}{\kilo \meter}$ we get a restriction on the time step of
\begin{equation*}
  \Delta t \leq \frac{\Delta x}{\lambda} \leq \SI{e-5}{\second},
\end{equation*}
so that $10^{5}$ time steps per second have to be taken.
For any realistic power network, computing so many time steps is at least difficult to implement and probably either too slow or too expensive to be practical.
On the other hand, the mean lifetime of the exponential decay of the Lyapunov function computed from~\eqref{eq:9}
is about $\SI{15}{\second}$ and this is a worst-case estimate as illustrated by~\ref{fig:lyapunov_net}.
Therefore for all practical purposes using the powerflow equations as is done in energy system simulation is a valid choice.

\section{Summary and future work}%
\label{sec:summary-future-work}
We have shown that the powerflow equations used in engineering applications can be thought of as special solutions to the Telegrapher's equations on the transmission lines between the nodes of the power network.
We also provided a second-order numerical scheme to compute these PDE solutions.
The scheme --- with possible exclusion of the coupling points--- is furthermore Lyapunov stable, which mimics the analytical solution. The Lyapunov stability 
including the coupling conditions will be investigated in future research using ideas presented in~\cite{kreiss_et_al_stability_boundary}.

Our numerical methods are extendable to higher order as there are splitting schemes at least up to order $6$\footnote{see \url{https://www.asc.tuwien.ac.at/~winfried/splitting/index.php}} and there are entropy-stable ENO schemes of arbitrary order~\cite{doi:10.1137/110836961} so that the Telegrapher's equation could be solved with much higher accuracy. However, comparing the computational costs, actually computing PDE solutions for realistic physical parameters can be only done for small instances and with high implementation effort.

\section*{Acknowledgments}
The authors gratefully thank the BMBF project ENets (05M18VMA) and the DFG project GO1920/10-1 for the financial support.

\printbibliography[heading=bibintoc]

\end{document}